# On the Exact Solution of Burgers-Huxley Equation Using the Homotopy Perturbation Method

S. Salman Nourazar[1], Mohsen Soori[1\*], Akbar Nazari-Golshan[2]

[1]Department of Mechanical Engineering, Amirkabir University of Technology (Tehran Polytechnic), Tehran, Iran
[2]Department of Physics, Amirkabir University of Technology (Tehran Polytechnic), Tehran, Iran
Email: [\*]mohsen.soori@gmail.com, [\*]m.soori@aut.ac.ir

## Abstract

The Homotopy Perturbation Method (HPM) is used to solve the Burgers-Huxley non-linear differential equations. Three case study problems of Burgers-Huxley are solved using the HPM and the exact solutions are obtained. The rapid convergence towards the exact solutions of HPM is numerically shown. Results show that the HPM is efficient method with acceptable accuracy to solve the Burgers-Huxley equation. Also, the results prove that the method is an efficient and powerful algorithm to construct the exact solution of non-linear differential equations.

## 1. Introduction

Most of the nonlinear differential equations do not have an analytical solution. Recently, semi-analytical solutions of real-life mathematical modeling are considered as a key tool to solve nonlinear differential equations.

The idea of the Homotopy Perturbation Method (HPM) which is a semi-analytical method was first pioneered by He [1]. Later, the method is applied by He [2] to solve the non-linear non-homogeneous partial differential equations. Nourazar *et al.* [3] used the homotopy perturbation method to find exact solution of Newell-Whitehead-Segel equation. Krisnangkura *et al.* [4] obtained exact traveling wave solutions of the generalized Burgers-

---

[\*]Corresponding author.

**How to cite this paper:** Nourazar, S.S., Soori, M. and Nazari-Golshan, A. (2015) On the Exact Solution of Burgers-Huxley Equation Using the Homotopy Perturbation Method. *Journal of Applied Mathematics and Physics*, **3**, 285-294.
http://dx.doi.org/10.4236/jamp.2015.33042

Huxley equation by using hyperbolic tangent method. New exact solutions of the generalized Burgers-Huxley equation are also obtained by Gao and Zhao [5]. Hashim *et al.* [6] have solved the generalized Burgers-Huxley equation by using adomian decomposition method. Wang *et al.* [7] studied the solitary wave solutions of the generalized Burgers-Huxley equation. Also, Darvishi *et al.* [8] have used spectral collocation method and Darvishi's preconditionings to solve the generalized Burgers-Huxley equation.

The generalized Burgers-Huxley equation is investigated by Satsuma [9] as:

$$\frac{\partial u}{\partial t} = \frac{\partial^2 u}{\partial x^2} - \alpha u^n \frac{\partial u}{\partial x} + \beta u(1-u^n)(u^n-\gamma), \quad 0 \leq x \leq 1, \ t \geq 0 \tag{1.1}$$

where $\alpha, \beta \geq 0$ are real constants and $n$ is a positive integer and $\gamma \in [0,1]$. Equation (1.1) models the interaction between reaction mechanisms, convection effects and diffusion transports [9].

When $\alpha = 0$ and $n = 1$, Equation (1.1) is reduced to the Huxley equation which describes nerve pulse propagationin nerve fibers and wall motion in liquid crystals [10]. When $\beta = 0$ and $n = 1$, Equation (1.1) is reduced to the Burgers equation describing the far field of wave propagation in nonlinear dissipative systems [11]. When $n = 1$, $\alpha \neq 0$ and $\beta \neq 0$, Equation (1.1) is turned into the Burgers-Huxley equation. The Burgers-Huxley equation is used to model the interaction between reaction mechanisms, convection effects and diffusion transport, nerve pulse propagation in nerve fibers as well as wall motion in liquid crystals.

In the present research work, the Homotopy Perturbation Method (HPM) is applied to obtain the closed form solution of the non-linear Burgers-Huxley equation. Three case study problems of non-linear Burgers-Huxley equations are solved by using the HPM. The trend of the rapid convergence towards the exact solution is shown when compared to the exact solution.

The idea of homotopy perturbation method is presented in Section 2. Application of the homotopy perturbation method to the exact solution of Burgers-Huxley equation is presented in Section 3.

## 2. The Idea of Homotopy Perturbation Method

The Homotopy Perturbation Method (HPM) is originally initiated by He [1]. This is a combination of the classical perturbation technique and homotopy technique. The basic idea of the HPM for solving nonlinear differential equations is as follow; consider the following differential equation:

$$E(u) = 0, \tag{2.1}$$

where $E$ is any differential operator. We construct a homotopy as follow:

$$H(u,p) = (1-p)F(u) + p(E(u) - F(u)). \tag{2.2}$$

where $F(u)$ is a functional operator with the known solution $v_0$. It is clear that when $p$ is equal to zero then $H(u,0) = F(u) = 0$, and when $p$ is equal to 1, then $H(u,1) = E(u) = 0$. It is worth noting that as the embedding parameter $p$ increases monotonically from zero to unity the zero order solution $v_0$ continuously deforms into the original problem $E(u) = 0$. The embedding parameter, $p \in [0,1]$, is considered as an expanding parameter [2]. In the homotopy perturbation method the embedding parameter $p$ is used to get series expansion for solution as:

$$u = \sum_{i=0}^{\infty} p^i v_i = v_0 + p v_1 + p^2 v_2 + p^3 v_3 + \cdots \tag{2.3}$$

when $p \to 1$, then Equation (2.2) becomes the approximate solution to Equation (2.1) as:

$$u = v_0 + v_1 + v_2 + v_3 + \cdots \tag{2.4}$$

The series Equation (2.4) is a convergent series and the rate of convergence depends on the nature of Equation (2.1) [1] [2]. It is also assumed that Equation (2.2) has a unique solution and by comparing the like powers of $p$ the solution of various orders is obtained. These solutions are obtained using the Maple package.

## 3. The Burgers-Huxley Equation

To illustrate the capability and reliability of the method, three cases of nonlinear diffusion equations are presented.

**Case I**: in this case we will examine the Burgers-Huxley equation for $\alpha = 0$, $n = 1$, $\gamma = 1$, $\beta = 1$, so, the equation is written as:

$$\frac{\partial u}{\partial t} = \frac{\partial^2 u}{\partial x^2} + u(1-u)(u-1) \tag{3.1}$$

Subject to initial condition:

$$u(x,0) = \frac{e^{\frac{\sqrt{2}}{4}x}}{e^{\frac{\sqrt{2}}{4}x} + e^{-\frac{\sqrt{2}}{4}x}} \tag{3.2}$$

We construct a homotopy for Equation (3.1) in the following form:

$$H(\upsilon, p) = (1-p)\left[\frac{\partial \upsilon}{\partial t} - \frac{\partial u_0}{\partial t}\right] + p\left[\frac{\partial \upsilon}{\partial t} - \frac{\partial^2 \upsilon}{\partial x^2} - \upsilon(1-\upsilon)(\upsilon-1)\right] \tag{3.3}$$

The solution of Equation (3.1) can be written as a power series in $p$ as:

$$\upsilon = \upsilon_0 + p\upsilon_1 + p^2\upsilon_2 + \cdots \tag{3.4}$$

Substituting Equation (3.4) and Equation (3.2) into Equation (3.3) and equating the terms with identical powers of $p$:

$$\begin{aligned}
&p^0: \frac{\partial \upsilon_0}{\partial t} = \frac{\partial u_0}{\partial t}, && \upsilon_0(x,0) = \frac{e^{\frac{\sqrt{2}}{4}x}}{e^{\frac{\sqrt{2}}{4}x} + e^{-\frac{\sqrt{2}}{4}x}}, \\
&p^1: \frac{\partial \upsilon_1}{\partial t} + \frac{\partial u_0}{\partial t} = \frac{\partial^2 \upsilon_0}{\partial x^2} + \upsilon_0(1-\upsilon_0)(\upsilon_0-1), && \upsilon_1(x,0) = 0, \\
&p^2: \frac{\partial \upsilon_2}{\partial t} = \frac{\partial^2 \upsilon_1}{\partial x^2} + 4\upsilon_0\upsilon_1 - \upsilon_1 - 3\upsilon_1\upsilon_0^2, && \upsilon_2(x,0) = 0, \\
&p^3: \frac{\partial \upsilon_3}{\partial t} = \frac{\partial^2 \upsilon_2}{\partial x^2} + 4\upsilon_0\upsilon_2 - \upsilon_2 - 3\upsilon_2\upsilon_0^2 + 2\upsilon_1^2 - 3\upsilon_0\upsilon_1^2, && \upsilon_3(x,0) = 0.
\end{aligned} \tag{3.5}$$

Using the Maple package to solve recursive sequences, Equation (3.5), we obtain the followings:

$$\upsilon_0(x,t) = \frac{e^{\frac{\sqrt{2}}{4}x}}{e^{\frac{\sqrt{2}}{4}x} + e^{-\frac{\sqrt{2}}{4}x}}, \quad \upsilon_1(x,t) = -\frac{1}{2}\frac{1}{\left(e^{\frac{\sqrt{2}}{4}x} + e^{-\frac{\sqrt{2}}{4}x}\right)^2}t,$$

$$\upsilon_2(x,t) = -\frac{1}{8}\frac{\left(e^{\frac{\sqrt{2}}{4}x} - e^{-\frac{\sqrt{2}}{4}x}\right)}{\left(e^{\frac{\sqrt{2}}{4}x} + e^{-\frac{\sqrt{2}}{4}x}\right)^3}t^2, \quad \upsilon_3(x,t) = -\frac{1}{48}\frac{\left(e^{\frac{\sqrt{2}}{4}x}\right)^2 - 4 + \left(e^{-\frac{\sqrt{2}}{4}x}\right)^2}{\left(e^{\frac{\sqrt{2}}{4}x} + e^{-\frac{\sqrt{2}}{4}x}\right)^4}t^3. \tag{3.6}$$

By setting $p = 1$ in Equation (3.4) the solution of Equation (3.1) can be obtained as $\upsilon = \upsilon_0 + \upsilon_1 + \upsilon_2 + \upsilon_3 + \cdots$
Therefore the solution of Equation (3.1) is written as:

$$\upsilon(x,t) = \frac{e^{\frac{\sqrt{2}}{4}x}}{e^{\frac{\sqrt{2}}{4}x} + e^{-\frac{\sqrt{2}}{4}x}} - \frac{1}{2}\frac{1}{\left(e^{\frac{\sqrt{2}}{4}x} + e^{-\frac{\sqrt{2}}{4}x}\right)^2}t - \frac{1}{8}\frac{\left(e^{\frac{\sqrt{2}}{4}x} - e^{-\frac{\sqrt{2}}{4}x}\right)}{\left(e^{\frac{\sqrt{2}}{4}x} + e^{-\frac{\sqrt{2}}{4}x}\right)^3}t^2 - \frac{1}{48}\frac{\left(e^{\frac{\sqrt{2}}{4}x}\right)^2 - 4 + \left(e^{-\frac{\sqrt{2}}{4}x}\right)^2}{\left(e^{\frac{\sqrt{2}}{4}x} + e^{-\frac{\sqrt{2}}{4}x}\right)^4}t^3 + \cdots \tag{3.7}$$

The Taylor series expansion for $\left(\dfrac{1}{2}+\dfrac{1}{2}\dfrac{e^{\frac{\sqrt{2}}{4}x-\frac{t}{4}}-e^{-\frac{\sqrt{2}}{4}x+\frac{t}{4}}}{e^{\frac{\sqrt{2}}{4}x-\frac{t}{4}}+e^{-\frac{\sqrt{2}}{4}x+\frac{t}{4}}}\right)$ is written as:

$$\frac{1}{2}+\frac{1}{2}\frac{e^{\frac{\sqrt{2}}{4}x-\frac{t}{4}}-e^{-\frac{\sqrt{2}}{4}x+\frac{t}{4}}}{e^{\frac{\sqrt{2}}{4}x-\frac{t}{4}}+e^{-\frac{\sqrt{2}}{4}x+\frac{t}{4}}}=\frac{e^{\frac{\sqrt{2}}{4}x}}{e^{\frac{\sqrt{2}}{4}x}+e^{-\frac{\sqrt{2}}{4}x}}-\frac{1}{2}\frac{1}{\left(e^{\frac{\sqrt{2}}{4}x}+e^{-\frac{\sqrt{2}}{4}x}\right)^{2}}t-\frac{1}{8}\frac{\left(e^{\frac{\sqrt{2}}{4}x}-e^{-\frac{\sqrt{2}}{4}x}\right)}{\left(e^{\frac{\sqrt{2}}{4}x}+e^{-\frac{\sqrt{2}}{4}x}\right)^{3}}t^{2}$$

$$-\frac{1}{48}\frac{\left(e^{\frac{\sqrt{2}}{4}x}\right)^{2}-4+\left(e^{-\frac{\sqrt{2}}{4}x}\right)^{2}}{\left(e^{\frac{\sqrt{2}}{4}x}+e^{-\frac{\sqrt{2}}{4}x}\right)^{4}}t^{3}+\cdots.$$

(3.8)

Combining Equation (3.8) with Equation (3.7), we get as follow:

$$\upsilon(x,t)=\frac{1}{2}+\frac{1}{2}\frac{e^{\frac{\sqrt{2}}{4}x-\frac{t}{4}}-e^{-\frac{\sqrt{2}}{4}x+\frac{t}{4}}}{e^{\frac{\sqrt{2}}{4}x-\frac{t}{4}}+e^{-\frac{\sqrt{2}}{4}x+\frac{t}{4}}}=\frac{1}{2}+\frac{1}{2}\tanh\left(\frac{1}{2\sqrt{2}}\left(x-\frac{t}{\sqrt{2}}\right)\right) \qquad (3.9)$$

This is the exact solution of the problem, Equation (3.1). **Table 1** shows the trend of rapid convergence of the results of $S_1(x,t)=v_0(x,t)$ to $S_6(x,t)=\sum_{i=0}^{5}v_i(x,t)$ using the HPM. The rapid convergence of the solution toward the exact solution, the maximum relative error of less than 0.0000058% is achieved as shown in **Table 1**.

**Case II**: In Equation (1.1) for $\alpha=-1$, $n=1$, $\gamma=1$, $\beta=1$, the Burgers-Huxley equation is written as:

$$\frac{\partial u}{\partial t}=\frac{\partial^{2}u}{\partial x^{2}}+u\frac{\partial u}{\partial x}+u(1-u)(u-1) \qquad (3.10)$$

**Table 1.** The percentage of relative errors of the results of $S_1(x,t)=v_0(x,t)$ to $S_6(x,t)=\sum_{i=0}^{5}v_i(x,t)$ of the HPM solution of Equation (3.1).

|  |  | Percentage of relative error (%RE) | | |
| --- | --- | --- | --- | --- |
|  |  | $x=1$ | $x=2$ | $x=3$ |
| $t=0.1$ | $S_1(x,t)$ | 0.01693168743 | 0.01002710463 | 0.005488150424 |
|  | $S_3(x,t)$ | 0.000002337346256 | 2.025644856e−7 | 9.527700575e−7 |
|  | $S_5(x,t)$ | 2.153484215e−10 | 3.464089104e−10 | 3.155246333e−10 |
|  | $S_6(x,t)$ | 9.744801271e−12 | 3.920421861e−11 | 9.937961084e−11 |
| $t=0.3$ | $S_1(x,t)$ | 0.05344388963 | 0.03164997413 | 0.01732302882 |
|  | $S_3(x,t)$ | 0.00006806597676 | 0.000003967422423 | 0.00002580410708 |
|  | $S_5(x,t)$ | 6.184344787e−8 | 9.467317545e−8 | 5.516785201e−8 |
|  | $S_6(x,t)$ | 1.008793018e−8 | 1.166422821e−9 | 2.026398250e−9 |
| $t=0.4$ | $S_1(x,t)$ | 0.07311570399 | 0.04329980772 | 0.02369935005 |
|  | $S_3(x,t)$ | 0.0001675410515 | 0.000007498069748 | 0.00006123279256 |
|  | $S_5(x,t)$ | 2.799248652e−7 | 4.009848425e−7 | 2.364525175e−7 |
|  | $S_6(x,t)$ | 5.775483086e−8 | 6.758498122e−9 | 1.111128458e−8 |

Subject to initial condition:

$$u(x,0) = \frac{e^{-\frac{x}{4}}}{e^{\frac{x}{4}} + e^{-\frac{x}{4}}} \tag{3.11}$$

To solve Equation (3.10) we construct a homotopy in the following form:

$$H(\upsilon, p) = (1-p)\left[\frac{\partial \upsilon}{\partial t} - \frac{\partial u_0}{\partial t}\right] + p\left[\frac{\partial \upsilon}{\partial t} - \frac{\partial^2 \upsilon}{\partial x^2} - \upsilon\frac{\partial \upsilon}{\partial x} - \upsilon(1-\upsilon)(\upsilon-1)\right] \tag{3.12}$$

The solution of Equation (3.10) can be written as a power series in $p$ as:

$$\upsilon = \upsilon_0 + p\upsilon_1 + p^2\upsilon_2 + \cdots \tag{3.13}$$

Substituting Equation (3.13) and Equation (3.11) in to Equation (3.12) and equating the term with identical powers of $p$, leads to:

$$\begin{aligned}
p^0 &: \frac{\partial \upsilon_0}{\partial t} = \frac{\partial u_0}{\partial t}, & \upsilon_0(x,0) &= \frac{e^{-\frac{x}{4}}}{e^{\frac{x}{4}} + e^{-\frac{x}{4}}}, \\
p^1 &: \frac{\partial \upsilon_1}{\partial t} + \frac{\partial u_0}{\partial t} = \frac{\partial^2 \upsilon_0}{\partial x^2} + \upsilon_0 \frac{\partial \upsilon_0}{\partial x} + \upsilon_0(1-\upsilon_0)(\upsilon_0 - 1), & \upsilon_1(x,0) &= 0, \\
p^2 &: \frac{\partial \upsilon_2}{\partial t} = \frac{\partial^2 \upsilon_1}{\partial x^2} + \upsilon_0 \frac{\partial \upsilon_1}{\partial x} + \upsilon_1 \frac{\partial \upsilon_0}{\partial x} + 4\upsilon_0\upsilon_1 - \upsilon_1 - 3\upsilon_1\upsilon_0^2, & \upsilon_2(x,0) &= 0, \\
p^3 &: \frac{\partial \upsilon_3}{\partial t} = \frac{\partial^2 \upsilon_2}{\partial x^2} + \upsilon_0 \frac{\partial \upsilon_2}{\partial x} + \upsilon_2 \frac{\partial \upsilon_0}{\partial x} + \upsilon_1 \frac{\partial \upsilon_1}{\partial x} + 4\upsilon_0\upsilon_2 - \upsilon_2 - 3\upsilon_2\upsilon_0^2 + 2\upsilon_1^2 - 3\upsilon_0\upsilon_1^2, & \upsilon_3(x,0) &= 0.
\end{aligned} \tag{3.14}$$

Using the Maple package to solve recursive sequences, Equation (3.14), we obtain the followings:

$$\upsilon_0(x,t) = \frac{e^{-\frac{x}{4}}}{e^{\frac{x}{4}} + e^{-\frac{x}{4}}}, \quad \upsilon_1(x,t) = -\frac{3}{4}\frac{1}{\left(e^{\frac{x}{4}} + e^{-\frac{x}{4}}\right)^2}t,$$

$$\upsilon_2(x,t) = \frac{9}{32}\frac{\left(e^{\frac{x}{4}} - e^{-\frac{x}{4}}\right)}{\left(e^{\frac{x}{4}} + e^{-\frac{x}{4}}\right)^3}t^2, \quad \upsilon_3(x,t) = -\frac{9}{128}\frac{\left(\left(e^{\frac{x}{4}}\right)^2 + \left(e^{-\frac{x}{4}}\right)^2 - 4\right)}{\left(e^{\frac{x}{4}} + e^{-\frac{x}{4}}\right)^4}t^3. \tag{3.15}$$

By setting $p = 1$ in Equation (3.13), the solution of Equation (3.10) can be obtained as $\upsilon = \upsilon_0 + \upsilon_1 + \upsilon_2 + \upsilon_3 + \cdots$
Therefore the solution of Equation (3.10) is written as:

$$\upsilon(x,t) = \frac{e^{-\frac{x}{4}}}{e^{\frac{x}{4}} + e^{-\frac{x}{4}}} - \frac{3}{4}\frac{1}{\left(e^{\frac{x}{4}} + e^{-\frac{x}{4}}\right)^2}t + \frac{9}{32}\frac{\left(e^{\frac{x}{4}} - e^{-\frac{x}{4}}\right)}{\left(e^{\frac{x}{4}} + e^{-\frac{x}{4}}\right)^3}t^2 - \frac{9}{128}\frac{\left(\left(e^{\frac{x}{4}}\right)^2 + \left(e^{-\frac{x}{4}}\right)^2 - 4\right)}{\left(e^{\frac{x}{4}} + e^{-\frac{x}{4}}\right)^4}t^3 + \cdots \tag{3.16}$$

The Taylor series expansion for $\left(\frac{1}{2} - \frac{1}{2}\frac{e^{\frac{x}{4}+\frac{3t}{8}} - e^{-\frac{x}{4}-\frac{3t}{8}}}{e^{\frac{x}{4}+\frac{3t}{8}} + e^{-\frac{x}{4}-\frac{3t}{8}}}\right)$ is written as:

$$\frac{1}{2}-\frac{1}{2}\frac{e^{\frac{x}{4}+\frac{3t}{8}}-e^{-\frac{x}{4}-\frac{3t}{8}}}{e^{\frac{x}{4}+\frac{3t}{8}}+e^{-\frac{x}{4}-\frac{3t}{8}}}=\frac{e^{-\frac{x}{4}}}{e^{\frac{x}{4}}+e^{-\frac{x}{4}}}-\frac{3}{4}\frac{1}{\left(e^{\frac{x}{4}}+e^{-\frac{x}{4}}\right)^2}t+\frac{9}{32}\frac{\left(e^{\frac{x}{4}}-e^{-\frac{x}{4}}\right)}{\left(e^{\frac{x}{4}}+e^{-\frac{x}{4}}\right)^3}t^2-\frac{9}{128}\frac{\left(\left(e^{\frac{x}{4}}\right)^2+\left(e^{-\frac{x}{4}}\right)^2-4\right)}{\left(e^{\frac{x}{4}}+e^{-\frac{x}{4}}\right)^4}t^3+\cdots \quad (3.17)$$

By substituting Equation (3.17) into Equation (3.16), Equation (3.16) can be reduced to:

$$\upsilon(x,t)=\frac{1}{2}-\frac{1}{2}\frac{e^{\frac{x}{4}+\frac{3t}{8}}-e^{-\frac{x}{4}-\frac{3t}{8}}}{e^{\frac{x}{4}+\frac{3t}{8}}+e^{-\frac{x}{4}-\frac{3t}{8}}}=\frac{1}{2}-\frac{1}{2}\tanh\left[\frac{1}{4}\left(x+\frac{3t}{2}\right)\right] \quad (3.18)$$

This is the exact solution of the problem, Equation (3.10). **Table 2** shows the trend of rapid convergence of the results of $S_1(x,t)=v_0(x,t)$ to $S_6(x,t)=\sum_{i=0}^{5}v_i(x,t)$ using the HPM solution toward the exact solution. The maximum relative error of less than 0.00014% is achieved in comparison to the exact solution as shown in **Table 2**.

**Case III**: In Equation (1.1) for $\alpha=-2$, $n=1$, $\gamma=3$, $\beta=1$, the Burgers-Huxley equation becomes:

$$\frac{\partial u}{\partial t}=\frac{\partial^2 u}{\partial x^2}+2u\frac{\partial u}{\partial x}+u(1-u)(u-3) \quad (3.19)$$

Subject to initial condition:

$$u(x,0)=\frac{3e^{\frac{-3(\sqrt{3}-1)}{4}x}}{e^{\frac{3(\sqrt{3}-1)}{4}x}+e^{-\frac{3(\sqrt{3}-1)}{4}x}} \quad (3.20)$$

**Table 2.** The percentage of relative errors of the results of $S_1(x,t)=v_0(x,t)$ to $S_6(x,t)=\sum_{i=0}^{5}v_i(x,t)$ of the HPM solution of Equation (3.10).

|  |  | Percentage of relative error (%RE) | | |
|---|---|---|---|---|
|  |  | $x=1$ | $x=2$ | $x=3$ |
| $t=0.1$ | $S_1(x,t)$ | 0.0484797171 | 0.056937877 | 0.063676094 |
|  | $S_3(x,t)$ | 0.0000184239461 | 0.000009125432 | 0.000006989146 |
|  | $S_5(x,t)$ | 7.8094040e−9 | 3.9049212e−9 | 1.37134980e−8 |
|  | $S_6(x,t)$ | 3.61301281e−10 | 6.1915442e−11 | 1.3095951e−10 |
| $t=0.3$ | $S_1(x,t)$ | 0.1570606291 | 0.184462686 | 0.206292613 |
|  | $S_3(x,t)$ | 0.000524561340 | 0.00023856284 | 0.00024584133 |
|  | $S_5(x,t)$ | 0.000000177771365 | 0.000000129074423 | 0.00000380612758 |
|  | $S_6(x,t)$ | 2.19608021e−7 | 2.10498615e−7 | 9.144252e−9 |
| $t=0.4$ | $S_1(x,t)$ | 0.2177728801 | 0.255767283 | 0.2860356311 |
|  | $S_3(x,t)$ | 0.001277016710 | 0.00055367038 | 0.00065831599 |
|  | $S_5(x,t)$ | 0.0000075558572 | 0.0000060480463 | 0.0000170599007 |
|  | $S_6(x,t)$ | 0.000001302473287 | 0.000001221861195 | 8.044206e−8 |

We construct a homotopy for Equation (3.19) in the following form:

$$H(v,p) = (1-p)\left[\frac{\partial v}{\partial t} - \frac{\partial u_0}{\partial t}\right] + p\left[\frac{\partial v}{\partial t} - \frac{\partial^2 v}{\partial x^2} - 2v\frac{\partial v}{\partial x} - v(1-v)(v-3)\right] \quad (3.21)$$

The solution of Equation (3.19) can be written as a power series in $p$ as:

$$v = v_0 + pv_1 + p^2 v_2 + \cdots \quad (3.22)$$

Substituting Equation (3.22) and Equation (3.20) into Equation (3.21) and equating the terms with identical powers of $p$:

$$p^0 : \frac{\partial v_0}{\partial t} = \frac{\partial u_0}{\partial t}, \quad v_0(x,0) = \frac{3e^{-\frac{3(\sqrt{3}-1)}{4}x}}{e^{\frac{3(\sqrt{3}-1)}{4}x} + e^{-\frac{3(\sqrt{3}-1)}{4}x}},$$

$$p^1 : \frac{\partial v_1}{\partial t} + \frac{\partial u_0}{\partial t} = \frac{\partial^2 v_0}{\partial x^2} + 2v_0\frac{\partial v_0}{\partial x} + v_0(1-v_0)(v_0 - 3), \quad v_1(x,0) = 0,$$

$$p^2 : \frac{\partial v_2}{\partial t} = \frac{\partial^2 v_1}{\partial x^2} + 2v_0\frac{\partial v_1}{\partial x} + 2v_1\frac{\partial v_0}{\partial x} + 8v_0 v_1 - 3v_1 - 3v_1 v_0^2, \quad v_2(x,0) = 0, \quad (3.23)$$

$$p^3 : \frac{\partial v_3}{\partial t} = \frac{\partial^2 v_2}{\partial x^2} + 2v_0\frac{\partial v_2}{\partial x} + 2v_2\frac{\partial v_0}{\partial x} + 2v_1\frac{\partial v_1}{\partial x} + 8v_0 v_2 - 3v_2$$
$$- 3v_2 v_0^2 + 4v_1^2 - 3v_0 v_1^2, \quad v_3(x,0) = 0.$$

Using the Maple package to solve recursive sequences, Equation (3.23), we obtain the followings:

$$v_0(x,t) = \frac{3e^{-\frac{3(\sqrt{3}-1)}{4}x}}{e^{\frac{3(\sqrt{3}-1)}{4}x} + e^{-\frac{3(\sqrt{3}-1)}{4}x}},$$

$$v_1(x,t) = -\frac{9}{2}\frac{(3\sqrt{3}-4)}{\left(e^{\frac{3(\sqrt{3}-1)}{4}x} + e^{-\frac{3(\sqrt{3}-1)}{4}x}\right)^2} t,$$

$$v_2(x,t) = \frac{27}{8}\frac{\left(e^{\frac{3(\sqrt{3}-1)}{4}x} - e^{-\frac{3(\sqrt{3}-1)}{4}x}\right)(43 - 24\sqrt{3})}{\left(e^{\frac{3(\sqrt{3}-1)}{4}x} + e^{-\frac{3(\sqrt{3}-1)}{4}x}\right)^3} t^2,$$

$$v_3(x,t) = \frac{27}{16}\left(\frac{1}{\left(e^{\frac{3(\sqrt{3}-1)}{4}x} + e^{-\frac{3(\sqrt{3}-1)}{4}x}\right)^4}\right)\left(\left(e^{\frac{3(\sqrt{3}-1)}{4}x}\right)^2 - 4 + \left(e^{-\frac{3(\sqrt{3}-1)}{4}x}\right)^2\right)(389 - 225\sqrt{3})t^3.$$

(3.24)

By setting $p = 1$ in Equation (3.22) the solution of Equation (3.19) can be obtained as $v = v_0 + v_1 + v_2 + v_3 + \cdots$
Thus the solution of Equation (3.19) can be written as:

$$v(x,t) = \frac{3e^{\frac{-3(\sqrt{3}-1)}{4}x}}{e^{\frac{3(\sqrt{3}-1)}{4}x} + e^{-\frac{3(\sqrt{3}-1)}{4}x}} - \frac{9}{2}\frac{(3\sqrt{3}-4)}{\left(e^{\frac{3(\sqrt{3}-1)}{4}x} + e^{-\frac{3(\sqrt{3}-1)}{4}x}\right)^2}t + \frac{27}{8}\frac{\left(e^{\frac{3(\sqrt{3}-1)}{4}x} - e^{-\frac{3(\sqrt{3}-1)}{4}x}\right)(43-24\sqrt{3})}{\left(e^{\frac{3(\sqrt{3}-1)}{4}x} + e^{-\frac{3(\sqrt{3}-1)}{4}x}\right)^3}t^2$$

(3.25)

$$+ \frac{27}{16}\left(\frac{1}{\left(e^{\frac{3(\sqrt{3}-1)}{4}x} + e^{-\frac{3(\sqrt{3}-1)}{4}x}\right)^4}\right)\left[\left(e^{\frac{3(\sqrt{3}-1)}{4}x}\right)^2 - 4 + \left(e^{-\frac{3(\sqrt{3}-1)}{4}x}\right)^2\right](389-225\sqrt{3})t^3 + \cdots.$$

The Taylor series expansion for $\left(\dfrac{3}{2} - \dfrac{3}{2}\dfrac{e^{\frac{3\sqrt{3}-3}{4}x + \frac{9\sqrt{3}-12}{4}t} - e^{-\frac{3\sqrt{3}-3}{4}x - \frac{9\sqrt{3}-12}{4}t}}{e^{\frac{3\sqrt{3}-3}{4}x + \frac{9\sqrt{3}-12}{4}t} + e^{-\frac{3\sqrt{3}+3}{4}x - \frac{9\sqrt{3}-12}{4}t}}\right)$ is written as:

$$\frac{3}{2} - \frac{3}{2}\frac{e^{\frac{3\sqrt{3}-3}{4}x+\frac{9\sqrt{3}-12}{4}t} - e^{-\frac{3\sqrt{3}-3}{4}x-\frac{9\sqrt{3}-12}{4}t}}{e^{\frac{3\sqrt{3}-3}{4}x+\frac{9\sqrt{3}-12}{4}t} + e^{-\frac{3\sqrt{3}+3}{4}x-\frac{9\sqrt{3}-12}{4}t}}$$

$$= \frac{3e^{\frac{-3(\sqrt{3}-1)}{4}x}}{e^{\frac{3(\sqrt{3}-1)}{4}x} + e^{-\frac{3(\sqrt{3}-1)}{4}x}} - \frac{9}{2}\frac{(3\sqrt{3}-4)}{\left(e^{\frac{3(\sqrt{3}-1)}{4}x} + e^{-\frac{3(\sqrt{3}-1)}{4}x}\right)^2}t + \frac{27}{8}\frac{\left(e^{\frac{3(\sqrt{3}-1)}{4}x} - e^{-\frac{3(\sqrt{3}-1)}{4}x}\right)(43-24\sqrt{3})}{\left(e^{\frac{3(\sqrt{3}-1)}{4}x} + e^{-\frac{3(\sqrt{3}-1)}{4}x}\right)^3}t^2 \quad (3.26)$$

$$+ \frac{27}{16}\left(\frac{1}{\left(e^{\frac{3(\sqrt{3}-1)}{4}x} + e^{-\frac{3(\sqrt{3}-1)}{4}x}\right)^4}\right)\left[\left(e^{\frac{3(\sqrt{3}-1)}{4}x}\right)^2 - 4 + \left(e^{-\frac{3(\sqrt{3}-1)}{4}x}\right)^2\right](389-225\sqrt{3})t^3 + \cdots.$$

Comparing Equation (3.26) with Equation (3.25), thus Equation (3.25) can be reduced to:

$$v(x,t) = \frac{3}{2} - \frac{3}{2}\frac{e^{\frac{3\sqrt{3}-3}{4}x+\frac{9\sqrt{3}-12}{4}t} - e^{-\frac{3\sqrt{3}-3}{4}x-\frac{9\sqrt{3}-12}{4}t}}{e^{\frac{3\sqrt{3}-3}{4}x+\frac{9\sqrt{3}-12}{4}t} + e^{-\frac{3\sqrt{3}+3}{4}x-\frac{9\sqrt{3}-12}{4}t}} = \frac{3}{2} - \frac{3}{2}\tanh\left[\frac{3\sqrt{3}-3}{4}\left(x + \frac{5-\sqrt{3}}{2}t\right)\right] \quad (3.27)$$

This is the exact solution of the problem, Equation (3.19). **Table 3** shows the trend of rapid convergence of the results of $S_1(x,t) = v_0(x,t)$ to $S_6(x,t) = \sum_{i=0}^{5} v_i(x,t)$ using the HPM solution toward the exact solution. The maximum relative error of less than 0.038% is achieved in comparison to the exact solution as shown in **Table 3**.

Deng [12] obtained some travelling solitary wave solutions of Equation (1.1) by applying the first-integral method as follows:

$$v(x,t) = \left[\frac{\gamma}{2} \pm \frac{\gamma}{2}\tanh\left(\frac{n\gamma(\rho \mp \alpha)}{4(n+1)}\left[x - \frac{(\alpha \mp \rho)\gamma + (\alpha \pm \rho)(n+1)}{2(n+1)}t + x_0\right]\right)\right]^{\frac{1}{n}} \quad (3.28)$$

S. S. Nourazar et al.

**Table 3.** The percentage of relative errors of the results of $S_1(x,t) = v_0(x,t)$ to $S_6(x,t) = \sum_{i=0}^{5} v_i(x,t)$ of the HPM solution of Equation (3.19).

|  |  | Percentage of relative error (%RE) | | |
|---|---|---|---|---|
|  |  | $x = 1$ | $x = 2$ | $x = 3$ |
| $t = 0.1$ | $S_1(x,t)$ | 0.1473751972 | 0.1768549738 | 0.1894968996 |
|  | $S_3(x,t)$ | 0.00008115001396 | 0.0004703355304 | 0.0008489173163 |
|  | $S_5(x,t)$ | 5.853207295e−7 | 0.000001136972415 | 2.175815927e−7 |
|  | $S_6(x,t)$ | 4.468762836e−8 | 5.878826669e−8 | 2.556253934e−8 |
| $t = 0.3$ | $S_1(x,t)$ | 0.5347070619 | 0.6416656548 | 0.6875331484 |
|  | $S_3(x,t)$ | 0.001445157793 | 0.01763206139 | 0.03053744832 |
|  | $S_5(x,t)$ | 0.0002129909216 | 0.0003484139717 | 0.00009299132649 |
|  | $S_6(x,t)$ | 0.00003726778757 | 0.00005691935052 | 0.00002679413014 |
| $t = 0.4$ | $S_1(x,t)$ | 0.7871664200 | 0.9446250035 | 1.012148624 |
|  | $S_3(x,t)$ | 0.002172311991 | 0.04914756514 | 0.08362907281 |
|  | $S_5(x,t)$ | 0.001086228788 | 0.001651183957 | 0.0005006633300 |
|  | $S_6(x,t)$ | 0.0002239213593 | 0.0003721014272 | 0.0001680489489 |

where

$$\rho = \sqrt{\alpha^2 + 4\beta(1+n)} \tag{3.29}$$

and $x_0$ is arbitrary constant.

This is in full agreement of the closed form solutions of Equation (3.1), Equation (3.10) and Equation (3.19) for differences value of parameters $\alpha$, $n$, $\gamma$ and $\beta$ in the three cases. So, it can be concluded that the HPM is a powerful and efficient technique to solve the non-linear Burgers-Huxley equation.

## 4. Conclusion

In the present research work, the exact solution of the Burgers-Huxley nonlinear diffusion equation is obtained using the HPM. The validity and effectiveness of the HPM is shown by solving three non-homogenous non-linear differential equations and the very rapid convergence to the exact solutions is also numerically demonstrated. The trend of rapid and monotonic convergence of the solution toward the exact solution is clearly shown by obtaining the relative error in comparison to the exact solution. The rapid convergence towards the exact solutions of HPM indicates that, using the HPM to solve the non-linear differential equations, a reasonable less amount of computational work with acceptable accuracy may be sufficient. Moreover, it can be concluded that the HPM is a very powerful and efficient technique which can construct the exact solution of nonlinear differential equations.

## References

[1] He, J.H. (1999) Homotopy Perturbation Technique. *Computer Methods in Applied Mechanics and Engineering*, **178**, 257-262. http://dx.doi.org/10.1016/S0045-7825(99)00018-3

[2] He, J.H. (2005) Application of Homotopy Perturbation Method to Nonlinear Wave Equations. *Chaos Solitons & Fractals*, **26**, 695-700. http://dx.doi.org/10.1016/j.chaos.2005.03.006

[3] Nourazar, S.S., Soori, M. and Nazari-Golshan, A. (2011) On the Exact Solution of Newell-Whitehead-Segel Equation

Using the Homotopy Perturbation Method. *Australian Journal of Basic and Applied Sciences*, **5**, 1400-1411.

[4]  Krisnangkura, M., Chinviriyasit, S. and Chinviriyasit, W. (2012) Analytic Study of the Generalized Burger's-Huxley Equation by Hyperbolic Tangent Method. *Applied Mathematics and Computation*, **218**, 10843-10847. http://dx.doi.org/10.1016/j.amc.2012.04.044

[5]  Gao, H. and Zhao, R.X. (2010) New Exact Solutions to the Generalized Burgers-Huxley Equation. *Applied Mathematics and Computation*, **217**, 1598-1603. http://dx.doi.org/10.1016/j.amc.2009.07.020

[6]  Hashim, I., Noorani, M.S.M. and Said Al-Hadidi, M.R. (2006) Solving the Generalized Burgers-Huxley Equation Using the Adomian Decomposition Method. *Mathematical and Computer Modelling*, **43**, 1404-1411. http://dx.doi.org/10.1016/j.mcm.2005.08.017

[7]  Wang, X.Y., Zhu, Z.S. and Lu, Y.K. (1990) Solitary Wave Solutions of the Generalised Burgers-Huxley Equation. *Journal of Physics A*: *Mathematical and General*, **23**, 271. http://dx.doi.org/10.1088/0305-4470/23/3/011

[8]  Darvishi, M.T., Kheybari, S. and Khani, F. (2008) Spectral Collocation Method and Darvishi's Preconditionings to Solve the Generalized Burgers-Huxley Equation. *Communications in Nonlinear Science and Numerical Simulation*, **13**, 2091-2103. http://dx.doi.org/10.1016/j.cnsns.2007.05.023

[9]  Satsuma, J. (1987) Topics in Soliton Theory and Exactly Solvable Nonlinear Equations. World Scientific, The Singapore City.

[10] Wang, X.Y. (1985) Nerve Propagation and Wall in Liquid Crystals. *Physics Letters A*, **112**, 402-406. http://dx.doi.org/10.1016/0375-9601(85)90411-6

[11] Whitham, G.B. (1974) Linear and Nonlinear Waves. Wiley, New York.

[12] Deng, X.J. (2008) Travelling Wave Solutions for the Generalized Burgers-Huxley Equation. *Applied Mathematics and Computation*, **204**, 733-737. http://dx.doi.org/10.1016/j.amc.2008.07.020